\documentstyle[12pt]{article}

\begin{document}

\def\ni{\noindent}

\def\vs{\vskip .3cm}

\def\ss{\smallskip}

\def\o{\over}

\def\la{\langle}

\def\ra{\rangle}

\def\.{\cdot}

\def\O{\Omega}

\def\d{\partial r}

\def\n{\nabla}

\def\l{\lambda}

\def\t{\tilde}

\def\be{\begin{equation}}

\def\ee{\end{equation}}

\def\beq{\begin{eqnarray*}}

\def\eeq{\end{eqnarray*}}



\parindent0em

\parskip1ex

\textwidth20cm



\newtheorem{de}{Definition}

\newtheorem{pr}{Proposition}

\newtheorem{th}[pr]{Theorem}

\newtheorem{lem}{Lemma}

\newtheorem{re}{Remark}

\newtheorem{cor}{Corollary}[pr]

\newtheorem{tcor}{Corollary}[th]

\newcommand{\f}{\mbox{f}}

\newcommand{\s}{\mbox{s} \,}

\newcommand{\rank}{\mbox{rank} \,}



{\Large \bf K\"ahlerian Killing Spinors, Complex Contact \\
 Structures and Twistor Spaces}\vskip 1cm

\ni A. Moroianu --{\footnotesize\sl Centre de Math\'ematiques, Ecole
Polytechnique, 91128 Palaiseau Cedex, France\\}\medskip

\ni U. Semmelmann --{\footnotesize\sl Humboldt-Universit\"at zu
Berlin, Institut f\"ur reine Mathematik (SFB 288), Ziegelstr. 13A,
D-10099 Berlin}\vskip 1cm

{\footnotesize

{\bf Abstract} - We collect our recent results ([5] and [8]) and we
get the equivalence of the three notions of the title under some
conditions. We then use this equivalence in order to prove some
consequences about Sasakian manifolds, complex almost contact
structures and complex k-contact structures.
}\vs

\vskip 1cm

\ni{\large \bf 0. Introduction.} The notion of a {\it complex contact
structure} was introduced in the late 50's by S. Kobayashi (cf. [6]),
in analogy to real contact structures.

In 1982 in [9], S. Salamon investigated quaternionic K\"ahler
manifolds. In particular, he defined the {\it twistor space} over
such a manifold as a generalization of the classical notion of
twistor space over a self-dual 4-manifold.

In 1986, K.D. Kirchberg was led to define {\it K\"ahlerian Killing
Spinors}, in order to characterize K\"ahler spin manifolds of odd
complex dimension admitting the smallest possible eigenvalue of the
Dirac operator (cf. [4]). Some important contributions to this
problem are also due to O. Hijazi (cf. [1]).

The aim of this paper is to collect our recent results (cf. [5] and
[8]), in order to explain the close connection between these  three
notions and to derive some corollaries. \vskip .5cm

\ni{\large \bf 1. Previous results.} In this section we describe the
three notions introduced above, and recall relevant results obtained
in each of these directions.

Let $M$ be a compact spin K\"ahler manifold of odd complex dimension
$m=n/2$ and positive scalar curvature $R$. Then, each eigenvalue
$\lambda$ of the Dirac operator $D$ satisfies the inequality (cf.
[4]) $$\lambda^2\ge{m+1\over 4m}\inf_M R.$$
In the limiting case of this inequality, $M$ is Einstein and any
eigenspinor $\Psi$ of $D$ corresponding to the eigenvalues
$\pm\sqrt{(m+1) R/4m}$ is a {\it K\"ahlerian Killing spinor},  i.e.,
satisfies the following first-order differential equation (cf. [1],
[4]):
$$\nabla_X\Psi+{1\over n+2}X\.D\Psi+{1\over n+2}J(X)\.\tilde D
\Psi=0.$$
We call such $M$ a {\it limiting manifold}. Conversely, any compact
K\"ahler manifold admitting K\"ahlerian Killing spinors is a limiting
manifold. The first known examples of such manifolds were the complex
projective spaces $CP^{2k+1}$.

Using complex contact structures it is possible to construct other
manifolds admitting K\"ahlerian Killing spinors. We will shortly
describe the
construction of [5].
\begin{de}{\em (cf.[6])} Let $ M^{2m} $ be a complex manifold of
complex dimension \\$ m = 2k+1 $. A {\it complex contact structure }
is a family ${\cal C } = \{ (U_i,\omega_i)\}$  satisfying the
following
conditions:
\begin{enumerate}
\renewcommand{\labelenumi}{{\em (\roman{enumi})} }
\item  $\{ U_i\}$  is an open covering of $ M $.
\item  $\omega_i$ is a holomorphic 1-form on $ U_i $.
\item $ \omega_i \wedge (\partial \omega_i)^k \in \Gamma (\Lambda
^{m,0}\,M )  $ is different from zero at every point of $U_i$.
\item $\omega_i = f_{ij} \omega_j $  in $ U_i \cap  U_j  \quad $,
where  $ f_{ij} $ is a holomorphic function on  $ U_i \cap  U_j .  $
\end{enumerate}
\end{de}
Let ${\cal C } = \{(U_i,\omega_i)\}$   be a complex contact
structure. Then there exists  an associated holomorphic line
subbundle $ L_{{\cal C }} \subset \Lambda^{1,0} (M) $ with transition
functions  $ \{ f_{ij}^{-1} \} $ and local sections  $\omega_i$. From
condition (iii)  immediately follows the isomorphism $ L_{{\cal C
}}^{k+1}  \; \cong \; K   ,  $ where $K = \Lambda^{m,0} (M)$ denotes
the canonical bundle of $M$. If we assume $k$ to be an odd integer
then $M$ admits a canonical spin  structure. It is given by the
isomorphism
\begin{equation}                                   \label{ident}
L _{{\cal C } }^{\frac{k+1}{2}}
\quad \cong \quad K^{\frac{1}{2}} \quad  \cong \quad  S_0  .
\end{equation}
Here $S_0$ is the subbundle of the spinor bundle $S$ which is defined
as  the eigenspace of  $\Omega$ for the eigenvalue $- i \,m $, where
the K\"ahler form $\Omega$ is considered as endomorphism of $S$. We
construct now a section $\Psi_{{\cal C } }$  of the spinor bundle
which is associated to the contact structure $\cal C$. For doing so
we fix $\; (U, \omega) \in {\cal C } \; $ and define  $\; \Psi_{{\cal
C } } \;$ over the open
set $U$ by
\begin{equation}        \label{spinor}
\Psi_{{\cal C } } \; \big|_{U }
:= | \Psi_{\omega} |^{-2} \; \bar{\eta}_{\omega} \; \cdot
\;\Psi_{\omega}  \;,
\end{equation}
where $ \Psi_{\omega} \in \Gamma ( S_0 \; \big|_{U }) $ is the local
section in $S_0$ corresponding to  $\omega^{\otimes \,
\frac{k+1}{2}}$ under the identification (\ref{ident}) and
$\eta_{\omega} := \omega \wedge (\partial\omega)^{\frac{k-1}{2}} $.
{}From the condition (iv)  it follows that the spinor $  \Psi_{{\cal
C
} }$ is globally defined. We have the following
\begin{pr}{\em (cf. [5])}                         \label{kks}
Let $ ( M,g,J) $ be a compact K\"ahler-Einstein manifold of complex
dimension $ m = 2k+1 $ with k odd, and let ${\cal C } $ be a complex
contact
structure on $ M $. Then the spinor $\Psi_{{\cal C } }$ associated
with   ${\cal C } $ satisfies the equation
$$ D^2 \, \Psi_{{\cal C } }  \; = \; \frac{m+1}{4m} \; R
\;\Psi_{{\cal C } }    ,$$
where $R$ is the scalar curvature of $(M,g)$. In particular, the
spinors  \\$ \Psi_{{\cal C } }^{\pm} := \lambda_1 \Psi_{{\cal C } }
\pm D \Psi_{{\cal C } } $ are  K\"ahlerian  Killing spinors, where $
\lambda_1  = \sqrt{\frac{m+1}{4m} R}$.
\end{pr}
A class of manifolds satisfying the assumptions of Proposition
\ref{kks} are the twistor spaces of quaternionic K\"ahler manifolds
introduced by S. Salamon (cf. [9]).

A {\it quaternionic K\"ahler manifold } is defined to be an oriented
4n-dimensional Riemannian manifold whose restricted holonomy group is
contained  in the subgroup $Sp(n) Sp(1) \subset SO(4n) \; (n \ge 2)$.
Salamon's idea is to construct over each such manifold $M$ a natural
$ C P^1$-- bundle  Z, admitting a K\"ahler metric such that the
bundle projection is a Riemannian submersion. He called this bundle
the  {\it twistor space } of $M$.
\begin{pr}{\em (cf.[9])}          \label{salamon}
Let  $ M^{4k} $ be a quaternionic K\"ahler manifold with positive
scalar curvature. Then its twistor space $Z$ admits a K\"ahler
Einstein metric of positive scalar curvature and a  complex contact
structure. Moreover, $Z$ is spin for odd k and $Z$ is spin for even k
iff $Z =  CP^{2k+1}$.
\end{pr}
{}From Propositions \ref{kks} and \ref{salamon} we obtain that all
the
twistor spaces of quaternionic K\"ahler manifolds $M^{4k}  \; (k
\equiv 1(2))$ with positive scalar curvature admits K\"ahlerian
Killing spinors, i.e. they are limiting manifolds.

The only explicitly known manifolds of this kind are the following
three families:
\begin{itemize}
\item $   \quad Sp(k+1) / Sp(k) \times  U(1)  \quad \cong \quad
CP^{2k+1}    $,
\item $   \quad SU(k+2) / S(U(k) \times  U(1) \times  U(1))     $,
\item $   \quad SO(k+4) / S( O(k) \times  O(3) \times  O(2) )   $.
\end{itemize}
and the 15--dimensional exceptional space  $ \quad F_4 / Sp(3) U(1)$.

It is now interesting to see that each such limiting manifold (i.e.
each spin K\"ahler manifold of odd complex dimension and positive
scalar curvature admitting K\"ahlerian Killing spinors) has to be a
twistor space. This is due to the following classification result:
\begin{pr}{\em (cf. [8])}       \label{andrei}
The limiting manifolds of complex dimension $4l+3$ are exactly the
twistor spaces associated to quaternionic K\"ahler manifolds of
positive scalar curvature. The only limiting manifold of complex
dimension $4l+1$ is $CP^{4l+1}$.
\end{pr}
The idea of the proof is the following. Take  a limiting manifold $M$
and consider  a maximal root of the canonical line bundle with some
hermitian metric. The associated principal $U(1)$-bundle over $M$,
say $P_M$, with a carefully chosen metric, is spin, and any spinor on
$M$ induces a {\it projectable} spinor on $P_M$. Moreover, a
K\"ahlerian Killing spinor induces a projectable real Killing spinor
on $P_M$. This forces $P_M$ to admits a regular Sasakian 3-structure
and $M$ to be the twistor space over the quotient of $P_M$ by the
Sasakian 3-structure.

The last part of the proposition follows from the fact that the only
spin twistor space of complex dimension 4l+1 is $CP^{4l+1}.$\vs

\ni{\large \bf 2. The results.} Combining the above propositions we
have
\begin{th}   \label{gesamt}
Let $M$ be a compact spin K\"ahler manifold of positive scalar
curvature and complex dimension $4l+3$. Then the following statements
are
equivalent:

{\em (i)} $\ \ M$ admits K\"ahlerian Killing spinors;

{\em (ii)} $\ M$ is K\"ahler-Einstein and admits a complex contact
structure;

{\em (iii)} $\ M$ is the twistor space of some quaternionic K\"ahler
manifold of positive scalar curvature.
\end{th}
As an immediate corollary we have the following result:
\begin{cor}                             \label{c1}
If $M$ is a K\"ahler-Einstein manifold of complex dimension $4l+3$
which admits a  complex contact structure, then $M$ is the twistor
space of some quaternionic K\"ahler manifold of positive scalar
curvature.
\end{cor}
This corollary is in fact part of a very recent theorem of C. LeBrun
(cf. [7]). He proves the same statement but without the restriction
on the dimension and also using a different method. The interest of
our proof lies in the unexpected appearance of the Dirac operator. As
a less obvious corollary we have the following
\begin{cor}                               \label{reg}
Let $M$ be a Riemannian manifold of real dimension $n=8l+7$,
admitting a Sasakian 3-structure which is regular in one direction.
Then it is regular in all directions.
\end{cor}
\ni{\sl Proof.} Let $V$ be the Killing vector field in the regular
direction. We denote by $N$ the quotient of $M$ by the $S^1$-action
in the direction of $V$. Regularity just means that $N$ is a
manifold. Now a simple calculation (cf. [2]) shows that $N$ is a
K\"ahler--Einstein manifold admitting a complex contact structure.

Corollary \ref{c1} yields that $N$ is the twistor space of some
quaternionic K\"ahler manifold $Q$, of positive scalar curvature.
Using [2] once again, we see that the 2-distribution given by the two
other Killing vector fields of the Sasakian 3-structure, projects on
the 2-distribution $\Theta$ which gives the complex contact structure
on $N$. So the quotient of $M$ by the Sasakian 3-structure is
diffeomorphic to the space of leaves of $\Theta$, which is exactly
the manifold $Q$. Thus our Sasakian 3-structure is regular.\vs
\begin{re}
Corollary {\em \ref{reg}} is also true for $n=8l+3$. We just have to
use the result of C. LeBrun {\em ([7])} instead of Corollary  {\em
\ref{c1}}
in the above proof.
\end{re}
In [3]  S. Ishihara and M. Konishi introduced the concept of {\it
complex almost contact structures}. These are the hermitian manifolds
of odd complex dimension $2n+1$ whose structure group can be reduced
to $U(1) \, \times \, (Sp(n) \, \otimes \, U(1))     $. They proved
that
each such manifold under an additional normality condition admits a
K\"ahler--Einstein metric and also a complex contact structure. In
[2] they also showed  the existence of a normal complex almost
contact structure on the  $S^1$--quotient of a 3--Sasakian space
which is  regular in one direction. From Theorem~\ref{gesamt} we then
have
\begin{cor}
Let $M$ be a complete Hermitian manifold with a complex almost
contact structure. Then the structure is normal iff $M$ is the
twistor space of some quaternionic K\"ahler manifold of positive
scalar curvature.
\end{cor}
To give a last application of Theorem \ref{gesamt} we consider a
generalization of complex contact structures. For this let ${\cal C}
= \{ U_i, \omega_i \}$  be a family of  (local) r--forms which again
satisfies conditions (i) -- (iv) of  Definition 1, where (iii)  has
to be changed  into:

{\em (iii)} ' $$ \quad \omega_i \wedge
(\partial \omega_i)^s  \in \Gamma (\Lambda ^{m,0}M  \, \big|_{U_i } )
\quad   \mbox{is different from zero at each point   of}
\quad U_i. $$
Here $s = \frac{m-r}{r+1}$ must be an integer. Such a family was
called a complex r--contact structure in [5]. If $s$ is an odd
integer then $M$ again admits a canonical spin  structure. In this
situation it is once more possible to construct a  K\"ahlerian
Killing spinor $\psi_{\cal C}$ (similar to  (\ref{spinor})).  Theorem
\ref{gesamt} then implies
\begin{pr}
Let $(M^{2m}, g, J)$ be a compact K\"ahler--Einstein manifold  with
positive scalar curvature which admits a complex r--contact structure
such
that $s = (m-r)/(r+1)$ is an odd integer. Then $M$ is a  complex
contact manifold.
\end{pr}
\begin{re}
As it was recently pointed out to us by K. Galicki, Corollary
\ref{reg} is an old result of S. Tanno [10].
\end{re}
\ni

{\footnotesize
{\sl Acknowledgements.} This note was finished during our stay at the
E. Schr\"odinger Institut in Vienna. We would like to thank the
institute for support and hospitality. }   \vs

\end{document}